\documentclass{article}
\usepackage[T1]{fontenc}
\usepackage{lmodern} 
\usepackage{amsmath,amsfonts,amssymb,amsthm}
\usepackage{graphicx}
\usepackage{xcolor}

\usepackage[letterpaper,top=2.7cm,bottom=2cm,left=2.8cm,right=2.8cm,marginparwidth=2cm]{geometry}

\usepackage[colorlinks=true, allcolors=blue]{hyperref}
\usepackage{caption}
\usepackage{subcaption}
\usepackage{booktabs}
\usepackage{float}
\usepackage{authblk}
\usepackage{threeparttable}
\usepackage{multirow}
\DeclareMathOperator{\thepace}{pace}

\newcommand{\LBpace}[1][]{\thepace^{LB}_{#1}}
\DeclareMathOperator{\npace}{npace}

\newcommand{\NLBpace}[1][]{\npace^{LB}_{#1}}
\newcommand{\init}{\texttt{root}}
\newcommand{\fin}{\texttt{end}}

\usepackage{natbib}
 \bibpunct[, ]{(}{)}{,}{a}{}{,}%
\title{Polynomial Optimization: \\ Enhancing RLT relaxations with Conic Constraints}
\author[1,2]{Brais Gonz\'alez-Rodr\'iguez\thanks{Corresponding Author: braisgonzalez.rodriguez@usc.es}}
\author[2]{Ra\'ul Alvite-Paz\'o}
\author[2]{Samuel Alvite-Paz\'o}
\author[3]{Bissan Ghaddar}
\author[1,2]{Julio Gonz\'alez-D\'iaz}
\affil[1]{Department of Statistics, Mathematical Analysis and Optimization and MODESTYA Research Group,  University of Santiago de Compostela.}
\affil[2]{CITMAga (Galician Center for Mathematical Research and Technology).}
\affil[3]{Ivey Business School, Western University, London, Ontario, Canada.}
\date{\today}

\begin{document}

\maketitle

\begin{abstract}
Conic optimization has recently emerged as a powerful tool for designing tractable and guaranteed algorithms for non-convex polynomial optimization problems. On the one hand, tractability is crucial for efficiently solving large-scale problems and, on the other hand, strong bounds are needed to ensure high quality solutions. In this research, we investigate the strengthening of RLT relaxations of polynomial optimization problems through the addition of nine different types of constraints that are based on linear, second-order cone, and semidefinite programming to solve to optimality the instances of well established test sets of polynomial optimization problems. We describe how to design these conic constraints and their performance with respect to each other and with respect to the standard RLT relaxations. Our first finding is that the different variants of nonlinear constraints (second-order cone and semidefinite) are the best performing ones in around $50\%$ of the instances. Additionally, we present a machine learning approach to decide on the most suitable constraints to add for a given instance. The computational results show that the machine learning approach significantly outperforms each and every one of the nine individual approaches.
\end{abstract}


\section{Introduction}

The large volume of theoretical and computational research on conic optimization has led to important advances over the last few years in the efficiency and robustness of the associated algorithmic procedures to solve them. Leading state-of-the-art mixed integer linear programming (MILP) solvers such as \texttt{Gurobi} \citep{gurobi}, \texttt{CPLEX} \citep{cplex}, and \texttt{Xpress} \citep{xpress} have recently added functionalities that allow to efficiently solve second-order cone programming (SOCP) problems and, further, \texttt{Mosek} \citep{mosekmanual,mosek} has positioned itself as a reliable solver for general semidefinite programming (SDP) problems.

Importantly, despite the convex nature of conic optimization problems, they are also proving to be a powerful tool in the design of branch-and-bound algorithms for general non-convex mixed integer nonlinear programming (MINLP) problems and, specially, polynomial optimization problems. The latter are very general and encompass many problems arising in operations research: integer programming, linear programming, mixed integer programming, and quadratic optimization, to name a few. The developments of the last several decades have shown that conic optimization is a central tool in addressing non-convexities. These advances have been more prominent in the case of polynomial optimization problems, and even more so in the particular case of (quadratically-constrained) quadratic optimization problems, where a variety of convex relaxations have been thoroughly studied \citep{shor1987, ghaddar2011, burer2020, bonami2019,elloumi2019}. These conic-based relaxations include semidefinite programming and second-order cone programming as their powerhouses. They often provide fast and guaranteed approaches for computing bounds on the global value of the non-convex optimization problem at hand. For instance,  \cite{lasserre2001} introduced semidefinite relaxations corresponding to liftings of the polynomial programs into higher dimensions. The construction is motivated by results related to representations of non-negative polynomials as sum-of-squares \citep{parrilo2003} and the dual theory of moments. This led to the development of systematic approaches for solving polynomial optimization problems to global optimality, the main limitation of these approaches currently being that they are computationally prohibitive in general.

Nowadays, state-of-the-art solvers tackle MINLP problems and, in particular, polynomial optimization ones through different relaxations to be solved at the nodes of the branch-and-bound tree. A crucial direction of current research focuses on integrating tighter and tighter relaxations while preserving reasonable computational properties, with relaxations that build upon different families of conic constraints becoming increasingly important. One relaxation for polynomial optimization problems is based on the Reformulation-Linearization Technique (RLT) and the most common approach to strengthen this relaxation comes from the identification of linear SDP-based cuts that tighten the RLT relaxation as the algorithm progresses \citep{Sherali2012cuts,Baltean2019,Gonzalez-rodriguez2020}. In this paper we are interested in the alternative approach of directly adding SOCP and SDP constraints to the relaxations as in \cite{burer2008} and \cite{buchheim2013}. Probably the main reason why this approach has received less attention so far is that the resulting algorithms need to rely on SDP solvers which, until recently, were probably not reliable enough and too expensive computationally in some instances. One of the primary goals of this paper is to show that the situation has changed, and that general branch-and-bound schemes based on the solution of SOCPs and SDPs at each and every node of the branch-and-bound tree can be very competitive and even superior to previous approaches.

In order to achieve the above goal, in this paper we introduce and compare linear SDP-based constraints as well as SOCP and SDP constraints. Importantly, their definition ensures that they are efficient in the sense of preserving the size and sparsity of the original RLT-based relaxation for polynomial optimization problems introduced in \cite{Sherali1992} and further refined in \cite{Sherali2013}. We consider a total of nine different versions of constraints to be added: one based on linear SDP-based cuts, four based on SOCP constraints, and four based on SDP constraints. These conic constraints are then integrated into the polynomial optimization solver \texttt{RAPOSa} \citep{Gonzalez-rodriguez2020}, whose core is an RLT-based branch-and-bound algorithm. As a second step, we then develop thorough computational studies on well established benchmarks including randomly generated instances from \cite{Dalkiran2016}, MINLPLib instances \citep{minlplib}, and QPLIB instances \citep{qplib}, which provide a wide variety of classes of polynomial optimization problems. One of the main findings is that, in around $50\%$ of the instances, the best performance is achieved by one of the versions explicitly incorporating SOCP and SDP conic constraints. The remaining $50\%$ is split quite evenly between the baseline RLT and the version that incorporates SDP-based linear cuts. Importantly, the analysis also allowed to identify particular classes of problems where one specific family of SOCP/SDP conic constraints is consistently superior to the linear versions. 

To the best of our knowledge, our contribution represents one of very few implementations of branch-and-bound schemes with conic relaxations for broad classes of problems, with the added value of the generality of the resulting scheme, since it can be applied to any given polynomial optimization problem. The most related approaches are \cite{burer2008} for non-convex quadratic problems with linear constraints and \cite{buchheim2013} for unconstrained mixed-integer quadratic problems.  In \cite{burer2008}, the authors present a computational analysis in which they compare the relative performance of different SDP relaxations, with the main highlight being that ``only a small number of nodes are required'' to fully solve the problems at hand. 
\cite{buchheim2013} introduce a new branch-and-bound algorithm, \texttt{Q-MIST} (Quadratic Mixed-Integer Semidefinite programming Technique), and develop a thorough computational analysis in which they show that, for a wide variety of families of instances within the scope of their algorithm, \texttt{Q-MIST} notably outperforms \texttt{Couenne} \citep{couenne}. It is worth noting that, if looking at specific classes of optimization problems, then one can find additional successful implementations of branch-and-bound algorithms with the inclusion of tailor-made conic constraints, such as \cite{ghaddar2019} for optimal power flow problems and combinatorial optimization problems, \cite{krislock2017} and \cite{rendl2010} for maximum cut problems, \cite{piccialli2022} for minimum sum-of-squares clustering, and \cite{ghaddar2011max} for maximum $k$-cut problems.

Last, but not least, in our computational experiments we also observe that there is a lot of variability in the best performing version of conic constraints for the different instances, with each version beating the rest for a non negligible number of instances. This observation motivates the last contribution of this paper, in which we exploit this variability by learning to choose the best version among our portfolio of different constraints. Building upon the framework in \cite{Bissan2022}, we show that the resulting machine learning version significantly outperforms each and every one of the underlying versions. This last contribution naturally fits into the rapidly emerging strand of research on ``learning to optimize'', whose advances are nicely presented in the survey papers \cite{Lodi2017} and \cite{Bengio2021}. The closest approach to this last part of our contribution is \cite{Baltean2019}, in which deep neural networks are used to rank SDP-based cuts for quadratic problems. Then, only the top scoring cuts are added, aiming to obtain a good balance between the tightness and the complexity of the relaxations. A common feature of the SDP-based cuts used in \cite{Baltean2019} and those in this paper is that they are designed with a big emphasis in sparsity considerations (number of nonzeros in the constraints), with the goal of obtaining computationally efficient relaxations. From the point of view of the learning process the approaches are quite different, since they focus on one type of constraints (the SDP-based cuts) and they want to learn the best SDP-based cuts to add at a given node, whereas we want learn to choose for any given instance which conic constraints to include in the relaxations.

The contribution of this paper can be summarized as follows. First, we define and show the potential of different SOCP and SDP strengthenings of the classic RLT relaxations for general polynomial optimization problems in a branch-and-bound scheme. Second, we design a machine learning approach to learn the best strenghtening to use on a given instance and obtain promising results for the resulting algorithm.

The remainder of this paper is organized as follows. In Section~\ref{sec:rlt} we present a brief overview of the classic RLT scheme. In Section~\ref{sec:conic} we describe the different families of conic constraints that will be integrated within the baseline RLT implementation. In Section~\ref{sec:results} we present a first series of computational results. Then, in Section~\ref{sec:learning} we show how the conic constraints can be further exploited within a machine learning framework. Finally, we conclude in Section~\ref{sec:conclusions} and discuss future research directions.


\section{Foundations of the RLT Technique}\label{sec:rlt}
The Reformulation-Linearization Technique was originally developed in \cite{Sherali1992}. It was designed to find global optima in polynomial optimization problems of the following form:
\begin{equation}
\begin{split}
\text{minimize} & \quad \phi_0(\mathbf{x})\\
\text{subject to}  & \quad \phi_r(\mathbf{x})\geq \beta_r, \quad r=1,2,\ldots, R_1 \\
& \quad \phi_r(\mathbf{x})=\beta_r, \quad r=R_1+1,\ldots,R\\
& \quad \mathbf{x}\in\Omega \subset \mathbb{R}^n\text{,}
\end{split}
\tag{\textbf{PO}}
\label{eq:PO}
\end{equation}
where $N = \lbrace 1, \dots, n \rbrace$ denotes the set of variables, each $\phi_r(\mathbf{x})$ is a polynomial of degree $\delta_r \in \mathbb{N}$ and $\Omega = \lbrace \mathbf{x} \in \mathbb{R}^n: 0 \leq l_j \leq x_j \leq u_j < \infty, \, \forall j \in N \rbrace \subset \mathbb{R}^n$ is a hyperrectangle containing the feasible region. Then, $\delta=\max_{r \in \{0,\ldots,R\}} \delta_r$ is the degree of the problem and $(N, \delta)$ represents all possible monomials of degree $\delta$.

The Reformulation-Linearization Technique consists of a branch and bound algorithm based on solving linear relaxations of the polynomial problem~\eqref{eq:PO}. These linear relaxations are built by working on a lifted space, where each monomial of the original problem is replaced with a corresponding RLT variable. For example, associated to monomials of the form $x_1x_2x_4$ and $x_1^2x_3^2$ one would define the RLT variables $X_{124}$ and $X_{1133}$, respectively. More generally, RLT variables are defined as
\begin{equation}
  X_J = \prod_{j \in J}x_j, 
  \label{eq:RLTidentity}
\end{equation}
where $J$ is a multiset containing the information about the multiplicity of each variable in the underlying monomial. Then, at each node of the branch-and-bound tree, one would solve the corresponding linear relaxation. Whenever we get a solution of a linear relaxation in which the identities in~\eqref{eq:RLTidentity} hold, we get a feasible solution of~\eqref{eq:PO}. Otherwise, the violations of these identities are used to choose the branching variable.

In order to get tighter relaxations and ensure convergence, new constraints, called bound-factor constraints, must be added. They are of the following form:
\begin{equation}\label{eq:boundfactor}
F_{\delta}(J_1,J_2)=\prod_{j\in J_1}{(x_j-l_j)}\prod_{j\in J_2}{(u_j-x_j)}\geq 0.
\end{equation}

Thus, for each pair of multisets $J_1$ and $J_2$ such that $J_1 \cup J_2 \subset (N,\delta)$ and $|J_1 \cup J_2| = \delta$, the corresponding bound-factor constraint is added to the linear relaxation.

\cite{Sherali2013} show that it is not necessary to add all bound-factor constraints to the linear relaxations, since certain subsets of them are enough to ensure the convergence of the algorithm. More precisely, they proved that convergence to a global optimum only requires the inclusion in the linear relaxation of those bound-factor constraints where $J_1 \cup J_2$ is a monomial that appears in \eqref{eq:PO}, regardless of its degree. Further, they also showed that convergence is also preserved if, whenever the bound-factor constraints associated to a monomial $J$ are present, all bound-factor constraints associated to monomials $J'\subset J$ are removed. Motivated by these results, $J$-sets are defined as those monomials of degree greater than one present in~\eqref{eq:PO} which, moreover, are not included in any other monomial (multiset inclusion). Consider the polynomial programming problem
\begin{equation}\label{eq:example}
	\begin{split}
	\text{minimize~~~~} & x_1^2 + x_2^2 + x_1x_2x_3 \\
	\text{subject to~~~~} & x_1x_2 + x_1x_4 \geq 1 \\
		& 1 \leq x_1 \leq 10 \\
		& 0 \leq x_2 \leq 8 \\
		& 0 \leq x_3 \leq 15 \\
		& 0 \leq x_4 \leq 7.
	\end{split}	
\end{equation}
The monomials with degree greater than one are $\{1,1\}$ ,$\{2,2\}$, $\{1,2,3\}$, $\{1,2\}$, and $\{1,4\}$. Since $\{1,2\}$ is included in $\{1,2,3\}$, it is removed. Therefore, the $J$-sets are $\{1,1\}$, $\{2,2\}$, $\{1,2,3\}$, and $\{1,4\}$. 

The analysis developed in this paper builds upon the above theoretical results from~\cite{Sherali2013} and, therefore, only the bound-factor constraints associated to $J$-sets are incorporated to the linear relaxations of~\eqref{eq:PO}. Such relaxations are less tight but, on the other hand, they are smaller in size and, hence, faster to solve. Note that the use of $J$-sets reduces not only the number of constraints but, more importantly, also the number of RLT variables. As it can be seen in~\cite{Sherali2013} and \cite{Gonzalez-rodriguez2020}, the performance of the RLT algorithm is clearly superior when the $J$-set approach is followed.

\section{Conic Enhancements of RLT}\label{sec:conic}

As already discussed earlier, the use of semidefinite programming to improve the performance of branch-and-bound schemes is not new, and \cite{Baltean2019} provides a thorough and up-to-date review of the field. Typically, the goal is to rely on semidefinite programming to tighten the relaxations of the original non-convex optimization problems targeted by a branch-and-bound algorithm. We start this section by reviewing the main ingredient of such methods and then present various families of SDP-driven constraints that can be incorporated into the RLT relaxations in an efficient way. In particular, they should preserve the underlying dimensionality and sparsity of the problem, which is crucial for these approaches to be competitive.

The main ingredient that has to be specified is the matrix or matrices on which positive semidefiniteness is to be imposed. To each (multi-)set of variables $\{x_j\}_{j\in J}$, with $J\subset (N,\delta)$, one can associate a vector $\omega=(x_j)_{j\in J}$ (resulting from the concatenation of all the variables in $J$, including repetitions). To any such vector one can associate matrix $M=\omega^T\omega$, which is trivially positive semidefinite. Now, let $ M_L = [\omega^T\omega]_L$ be the matrix obtained when each monomial in $M$ is replaced by the corresponding RLT variable in the lifted space. The constraint $M_L\succcurlyeq 0$ is a valid cut because it never removes feasible solutions of~\eqref{eq:PO} and, hence, it does not compromise convergence of the RLT algorithm to a global optimum. In practice, vectors of the form $(1,(x_j)_{j\in J})$ are often preferred, since they result in $M_L$ matrices containing also variables in the original space and not only RLT variables, leading to tighter relaxations. In \cite{Sherali2012cuts}, for instance, the authors discuss different ways of defining $w$ for a given $J$. Specifically, they mainly work with $\omega^1=(x_j)_{j\in J}$, $\omega^2=(1,(x_j)_{j\in J})$, and $\omega^3=(1, x_1, x_2, \ldots)$, where $\omega^3$ is defined by concatenating also monomials of degree greater than one, while ensuring that no monomial in the resulting matrix $M$ has degree larger than $\delta$, the degree of~\eqref{eq:PO}.

We now move to the definition of the specific SDP-driven constraints for the RLT algorithm that constitute the subject of study in this paper.

\subsection{Linear SDP-based Constraints}\label{sec:sdpcuts}
In \cite{Sherali2012cuts}, the authors associate linear cuts to the constraints of the form~$M_L\succcurlyeq 0$ as follows. At each node of the branch-and-bound tree, the positive semidefiniteness of the chosen $M_L$ matrices is assessed at the solution of the corresponding relaxation. Given a negative eigenvalue of one such matrix with $\alpha$ as its associated eigenvector, then the valid cut $\alpha^T M_L\alpha\geq 0$ can be added to the linear relaxation to separate the current solution. A potential drawback of these cuts is that they may be very ``dense'', in the sense of involving a large number of variables, which may increase the solving time of the relaxations. Thus, as already discussed in \cite{Sherali2012cuts}, it is important to carefully choose $M_L$ matrices.

The sparsity of the cuts is particularly important if the RLT algorithm is being run with the $J$-set approach since, in general, the resulting cuts might involve monomials not contained in any $J$-set. This would require to include additional RLT variables in the relaxations (and the corresponding bound-factor constraints), increasing the size and potentially reducing the sparsity of the relaxations. Here we follow \cite{Gonzalez-rodriguez2020}, where the authors consider, at each node, all $M_L$ matrices obtained from vectors $\omega^k$ associated with the different $J$-sets of~\eqref{eq:PO}, which lead to sparse cuts that essentially preserve the dimensionality of the resulting relaxations. The authors present a detailed computational analysis, comparing different approaches to add an inherit cuts. We adopt the best performing version, which consists of using vector $\omega^2$ and inheriting all cuts from one node to all its descendants.

In \cite{Sherali2012cuts} different methods are discussed to efficiently look for the negative eigenvalues of the $M_L$ matrices. One such approach, that we follow here, consists of dividing each matrix in $10\times 10$ overlapping submatrices (each matrix shares its first 5 rows with the preceding one) and add a valid cut for each negative eigenvalue they have. This approach, on top of being computationally cheap, leads to even sparser cuts.

The above discussion regarding the adequacy of building constraints based on $J$-sets, in order to obtain sparser constraints (number of nonzero coefficients) and to preserve the size (number of variables) of the resulting relaxations, also applies to the conic constraints defined in the following subsections which, therefore, also build upon $J$-sets.

\subsection{SDP Constraints}\label{sec:sdp}

We next describe two approaches to tighten the classic RLT relaxations by directly adding semidefinite constraints. They just differ in the matrices on which positive semidefiniteness is imposed.

\begin{description}
\item[Approach 1.] For each $J$-set $J$, $\omega^1$ is used to define the $M_L$ matrix and the constraint $M_L\succcurlyeq 0$. Thus, $M = (\omega^1)^T\omega^1$. Note that, whenever $J$ contains only one variable $x_i$ (possibly multiple times), this would result in a trivial constraint and, therefore, these constraints are disregarded with one exception: if $|J|=2$, then $\omega^1$ is replaced with $(1,x_i)$. With this exception, this approach is mathematically equivalent for quadratic problems to the SOCP approach we present in Section~\ref{sec:socc} below.
\item[Approach 2.] For each $J$-set $J$,  $\omega^2$ is used to define the $M_L$ matrix and the constraint $M_L\succcurlyeq 0$.
\end{description}

Preliminary analysis have shown that constraints building upon $\omega^3$, although they lead to tighter relaxations, generate significantly bigger matrices and increase the complexity of the resulting SDP problems.
To provide an example of both approaches, consider the polynomial optimization problem in Equation~\eqref{eq:example}. Then, the semidefinite constraints added with the above approaches are the following ones:

\begin{center}
{\footnotesize
\begin{tabular}{c}
     \multicolumn{1}{l}{\textbf{\normalsize Approach 1}}  \\
     $\displaystyle \begin{pmatrix}
    1 & x_1\\
   x_1 & X_{11}\\
    \end{pmatrix}
 \succcurlyeq  0, ~
 \begin{pmatrix}
    1 & x_2\\
   x_2 & X_{22}\\
    \end{pmatrix}
 \succcurlyeq  0, ~
 \begin{pmatrix}
    X_{11} & X_{12} & X_{13}\\
   X_{12} & X_{22} & X_{23}\\
   X_{13} & X_{23} & X_{33}\\
    \end{pmatrix}
 \succcurlyeq  0, \text{ and } ~
 \begin{pmatrix}
    X_{11} & X_{14}\\
   X_{14} & X_{44}\\
    \end{pmatrix}
 \succcurlyeq  0$.\\
    \multicolumn{1}{l}{\textbf{\normalsize Approach 2}}  \\
$\displaystyle \begin{pmatrix}
    1 & x_1 & x_1\\
   x_1 & X_{11} & X_{11}\\
   x_1 & X_{11} & X_{11}\\
    \end{pmatrix}
    \succcurlyeq  0, ~
    \begin{pmatrix}
    1 & x_2 & x_2\\
   x_2 & X_{22} & X_{22}\\
   x_2 & X_{22} & X_{22}\\
    \end{pmatrix}
 \succcurlyeq  0, ~
 \begin{pmatrix}
    1 & x_1 & x_2 & x_3\\
   x_1 & X_{11} & X_{12} & X_{13}\\
   x_2 & X_{12} & X_{22} & X_{23}\\
   x_3 & X_{13} & X_{23} & X_{33}\\
    \end{pmatrix}
 \succcurlyeq  0,  \text{ and } ~
 \begin{pmatrix}
    1 & x_1 & x_4\\
   x_1 & X_{11} & X_{14}\\
   x_4 & X_{14} & X_{44}\\
    \end{pmatrix}
 \succcurlyeq  0$.
\end{tabular}}
\end{center}


\subsection{SOCP Constraints}\label{sec:socc}
We now describe the second-order cone constraints which, with respect to the SDP ones, lead to looser relaxations but, on the other hand, can be solved more efficiently by state-of-the-art optimization solvers. For each $J$-set $J$ and each pair of variables present in $J$, $x_i\neq x_j$, we define the following second-order cone constraint:
\begin{align}\label{eq:socc1}
\frac{X_{ii} + X_{jj}}{2} &\geq \Biggl\lVert \begin{pmatrix}
    X_{ij}\\
    \displaystyle\frac{X_{ii} - X_{jj}}{2}\\
    \end{pmatrix} \Biggr\rVert_2.
\end{align}
We argue now why these constraints are valid cuts, \emph{i.e.,} they never remove solutions feasible to~\eqref{eq:PO}. Constraint~\eqref{eq:socc1} can be equivalently rewritten as $X_{ii}X_{jj} \geq X_{ij}^2$. Then, given a solution of a linear relaxation satisfying the RLT identities in~\eqref{eq:RLTidentity}, the above condition reduces to $x_ix_ix_jx_j\geq x_ix_jx_ix_j$, which is trivially true. Note that constraints in~\eqref{eq:socc1} are trivially true if $i=j$ and, hence, whenever we have a variable $x_i$ appearing twice or more in $J$, we instead add the second-order constraint
\begin{align}\label{eq:socc2}
\frac{1 + X_{ii}}{2} &\geq \Biggl\lVert \begin{pmatrix}
    x_i\\
    \displaystyle\frac{1 - X_{ii}}{2}\\
    \end{pmatrix} \Biggr\rVert_2,
\end{align}
which is equivalent to $X_{ii}\geq x_ix_i$ and, for solutions satisfying~\eqref{eq:RLTidentity}, is again trivially true. Consider again the polynomial optimization problem in Equation~\eqref{eq:example}. Then, the SOCP constraints added are the following ones:

{\footnotesize
\begin{equation*}
\frac{1 + X_{11}}{2} \geq \Biggl\lVert \begin{pmatrix}
    x_1\\
    \displaystyle\frac{1 - X_{11}}{2}\\
    \end{pmatrix} \Biggr\rVert_2, ~
\frac{1 + X_{22}}{2} \geq \Biggl\lVert \begin{pmatrix}
    x_2\\
    \displaystyle\frac{1 - X_{22}}{2}\\
    \end{pmatrix} \Biggr\rVert_2, ~
\frac{X_{11} + X_{22}}{2} \geq \Biggl\lVert \begin{pmatrix}
    X_{12}\\
    \displaystyle\frac{X_{11} - X_{22}}{2}\\
    \end{pmatrix} \Biggr\rVert_2, 
\end{equation*}
\begin{equation*}
\frac{X_{11} + X_{33}}{2} \geq \Biggl\lVert \begin{pmatrix}
    X_{13}\\
    \displaystyle\frac{X_{11} - X_{33}}{2}\\
    \end{pmatrix} \Biggr\rVert_2, ~
\frac{X_{22} + X_{33}}{2} \geq \Biggl\lVert \begin{pmatrix}
    X_{23}\\
    \displaystyle\frac{X_{11} - X_{33}}{2}\\
    \end{pmatrix} \Biggr\rVert_2, \text{ and } \,
\frac{X_{11} + X_{44}}{2} \geq \Biggl\lVert \begin{pmatrix}
    X_{14}\\
    \displaystyle\frac{X_{11} - X_{44}}{2}\\
    \end{pmatrix} \Biggr\rVert_2.
\end{equation*}
}

\subsection{Binding SOCP and SDP Constraints}

Since solving SOCP or SDP problems is usually more time-consuming than solving linear programming problems, we define a new approach in order to reduce the time needed for solving the resulting RLT relaxation with SOCP or SDP constraints. This consists of checking which conic constraints (second-order cone or semidefinite) are binding after solving the first relaxation, \emph{i.e.}, this is done only once, at the root node. Thereafter, only these binding constraints are used to tighten the future linear relaxations. This approach significantly reduces the number of second-order cone or semidefinite constraints in the relaxations and, although these new relaxations are not as tight, one might expect that the binding constraints at the root node tend to be the most important ones in subsequent relaxations, at least in the first phase of the algorithm. We assess the trade-off between the difficulty of solving the relaxations and how tight they are in the computational analysis in the next section.

\section{Computational Results} \label{sec:results}

\subsection{Testing Environment}
All the computational analyses reported in this paper have been performed on the supercomputer Finisterrae~III, provided by Galicia Supercomputing Centre (CESGA). Specifically, we use nodes powered with 32 cores Intel Xeon Ice Lake 8352Y CPUs with 256GB of RAM connected through an Infiniband HDR network, and 1TB of SSD.

Regarding the datasets, we use three different sets of problems. The first one, DS, is taken from \cite{Dalkiran2016} and consists of 180 instances of randomly generated polynomial programming problems of different degrees, number of variables, and density. The second dataset comes from the well known benchmark MINLPLib \citep{minlplib}, a library of Mixed-Integer Nonlinear Programming problems. We have selected from  MINLPLib those instances that are polynomial programming problems with box-constrained and continuous variables, resulting in a total of 166 instances. The third dataset comes from another well known benchmark, QPLIB \citep{qplib}, a library of quadratic programming instances, for which we made a selection analogous to the one made for MINLPLib, resulting in a total of 63 instances. Hereafter we refer to the first dataset as DS, to the second one as MINLPLib, and to the third one as QPLIB.\footnote{Instances from DS dataset can be downloaded at \url{https://raposa.usc.es/files/DS-TS.zip}, instances from MINLPLib dataset can be downloaded at \url{https://raposa.usc.es/files/MINLPLib-TS.zip}, and instances from QPLIB dataset can be downloaded at \url{https://raposa.usc.es/files/QPLIB-TS.zip}.}

We develop our analysis by building upon the global solver for polynomial optimization problems \texttt{RAPOSa} \citep{Gonzalez-rodriguez2020}. Regarding the auxiliary solvers, \texttt{RAPOSa} uses i)~\texttt{Gurobi} for the linear relaxations ii)~\texttt{Gurobi} or \texttt{Mosek} for the SOCP relaxations, and iii)~\texttt{Mosek} for the SDP relaxations. The main objective of the thorough numerical analysis developed in this and in the following section is to assess the performance of different SOCP/SDP conic-driven versions of \texttt{RAPOSa} with respect to two more traditional ones: basic RLT and RLT with linear SDP-based cuts. More precisely, the full set of ten different versions is as follows: 
\begin{itemize}
    \item \texttt{RLT}: standard RLT algorithm (with $J$-sets).
    \item \texttt{SDP-Cuts}: linear SDP-based cuts added to RLT.
    \item \texttt{SOCP$^G$}: SOCP constraints added to the RLT relaxation and solved with \texttt{Gurobi}.
    \item \texttt{SOCP$^{G,B}$}: same as above, but using only constraints that were binding at the root node.
    \item \texttt{SOCP$^M$}: SOCP constraints added to the RLT relaxation and solved with \texttt{Mosek}.
    \item \texttt{SOCP$^{M,B}$:} same as above, but using only constraints that were binding at the root node.
    \item \texttt{SDP$^1$}: SDP constraints added to the RLT following Approach~1.
    \item \texttt{SDP$^{1,B}$}: same as above, but using only constraints that were binding at the root node.
    \item \texttt{SDP$^{2}$}: SDP constraints added to the RLT following Approach~2.
    \item \texttt{SDP$^{2,B}$}: same as above, but using only constraints that were binding at the root node.
\end{itemize}

For each instance and each one of the above versions, we run \texttt{RAPOSa} with a time limit of one hour.

\subsection{Numerical Results and Analysis}

The main goal of the numerical analysis in this section is to show the potential of SOCP/SDP conic-constraints to improve upon the performance of more classic implementations of RLT, such as \texttt{RLT} and \texttt{SDP-Cuts}. The measures used to evaluate the performance of the different versions of \texttt{RAPOSa} are $\LBpace$ and $\NLBpace$, two performance indicators introduced in \cite{Bissan2022} that capture the pace at which a given algorithm closes the gap or, more precisely, the pace at which it increases the lower bound along the branch-and-bound tree. To compute $\LBpace$ we use the following formula:
\begin{equation}
    \LBpace=\frac{\text{time}}{LB^\fin-LB^\init+\varepsilon}.
\label{eq:epsilon}
\end{equation}
Then, $\NLBpace$ is just a normalized version of $\LBpace$ with values in $[0,1]$. It is computed, for each version of the solver/algorithm, by dividing the best (smallest) pace among all versions to be compared by the pace of the current one.

As thoroughly discussed in \cite{Bissan2022}, $\LBpace$ and $\NLBpace$ are natural measures that allow to compare the performance of different solvers/algorithms on all the instances of a test set at once, regardless of their difficulty and of how many versions of the underlying solver/algorithm have solved them to optimality. This is different from more common approaches, where the running time is used to evaluate performance on instances solved by all versions, the optimality gap is used for those instances solved by none, and where some decision has to be made regarding those instances solved by some but not all of the versions of the solver/algorithm.

\begin{figure}[!htbp]
\centering
\includegraphics[width=0.75\textwidth]{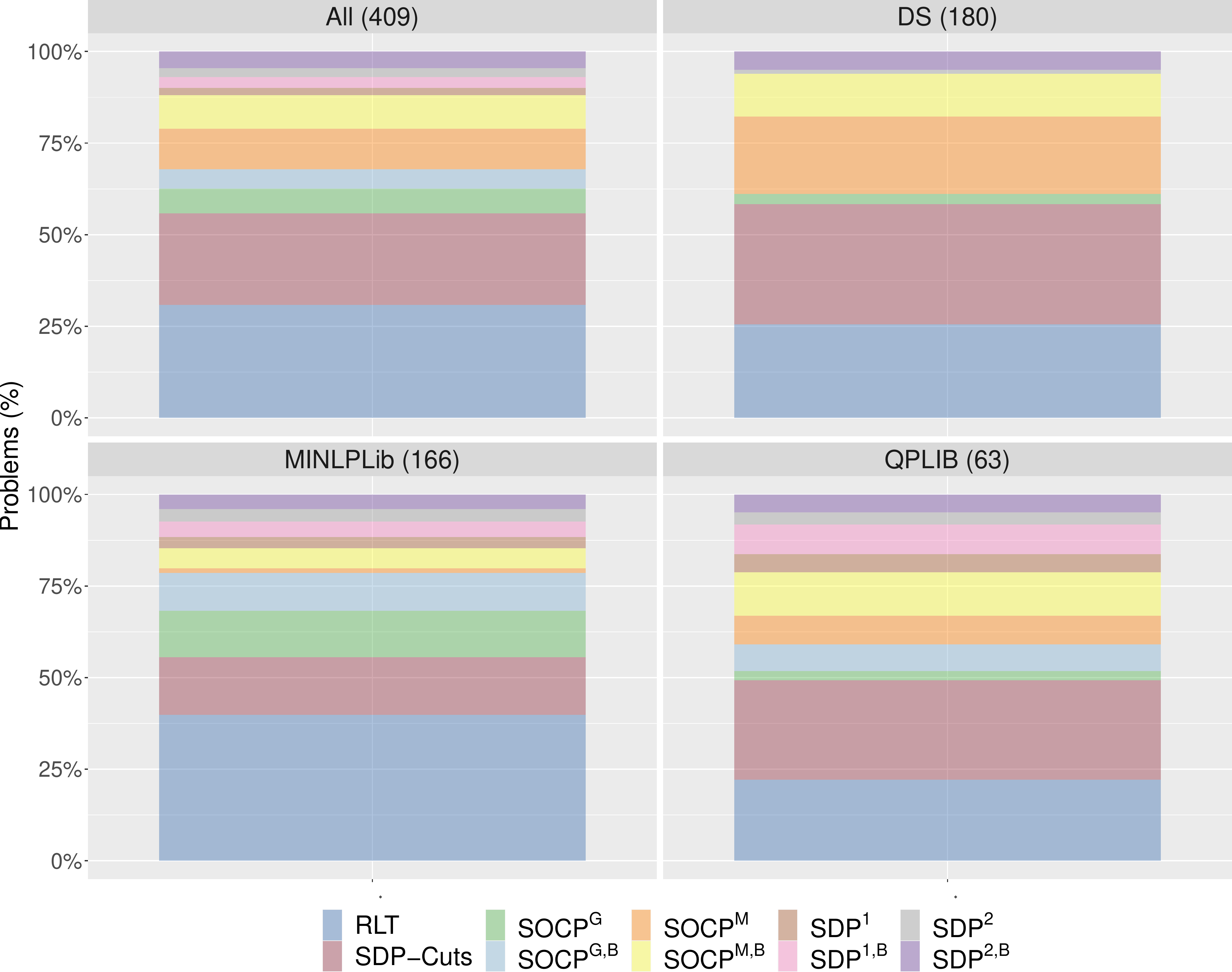}
\caption{Percentage of instances in which each version delivers the best performance.}
\label{fig:bar_comp_all_bat}
\end{figure}

Figure~\ref{fig:bar_comp_all_bat} shows, for the different sets of problems, the percentage of instances in which each one of the ten versions is the best one. We can see that, quite consistently across the three test sets, a version with either SOCP or SDP constraints is the best one in around $50\%$ of the instances. This is in itself one of the main highlights of this paper: RLT versions incorporating nonlinear SOCP/SDP conic constraints can improve the performance of RLT-based algorithms in half of the instances of the sets of problems under consideration. 

Figure~\ref{fig:bar_comp_all_bat} contains more valuable information. First, the versions with SOCP constraints are the best ones much more often than those with SDP constraints. Regarding the solvers for SOCP versions, \texttt{Gurobi} performs notably better in MINLPLib instances, whereas \texttt{Mosek} is overwhelmingly better on DS and is also the best one for QPLIB instances. When comparing binding versions with their non-binding counterparts we can see that, for SDP versions, \texttt{SDP$^{1,B}$} and \texttt{SDP$^{2,B}$} are the best ones significantly more often than \texttt{SDP$^{1}$} and \texttt{SDP$^{2}$}, respectively. In the case of SOCP versions, there is no clear winner between binding and non-binding versions. Finally, regarding the RLT versions without conic constraints, \texttt{RLT} and \texttt{SDP-Cuts}, we can see that each of them turns out to be the best option in around $25\%$ of the instances, with \texttt{SDP-Cuts} looking preferable in DS and QPLIB, whereas \texttt{RLT} is the best one three times as much in MINLPLib. 

Importantly, Figure~\ref{fig:bar_comp_all_bat} and the preceding discussion show that there is a lot of variability, with all ten versions showing up as the best choice for a non-negligible percentage of instances. Further, this variability also follows different patterns for the different sets of problems, which motivates the approach taken in Section~\ref{sec:learning} below, where we use machine learning techniques to try to learn to choose in advance the most promising RLT version for a given instance. 

\begin{figure}[!htbp]
\centering
\includegraphics[width=0.85\textwidth]{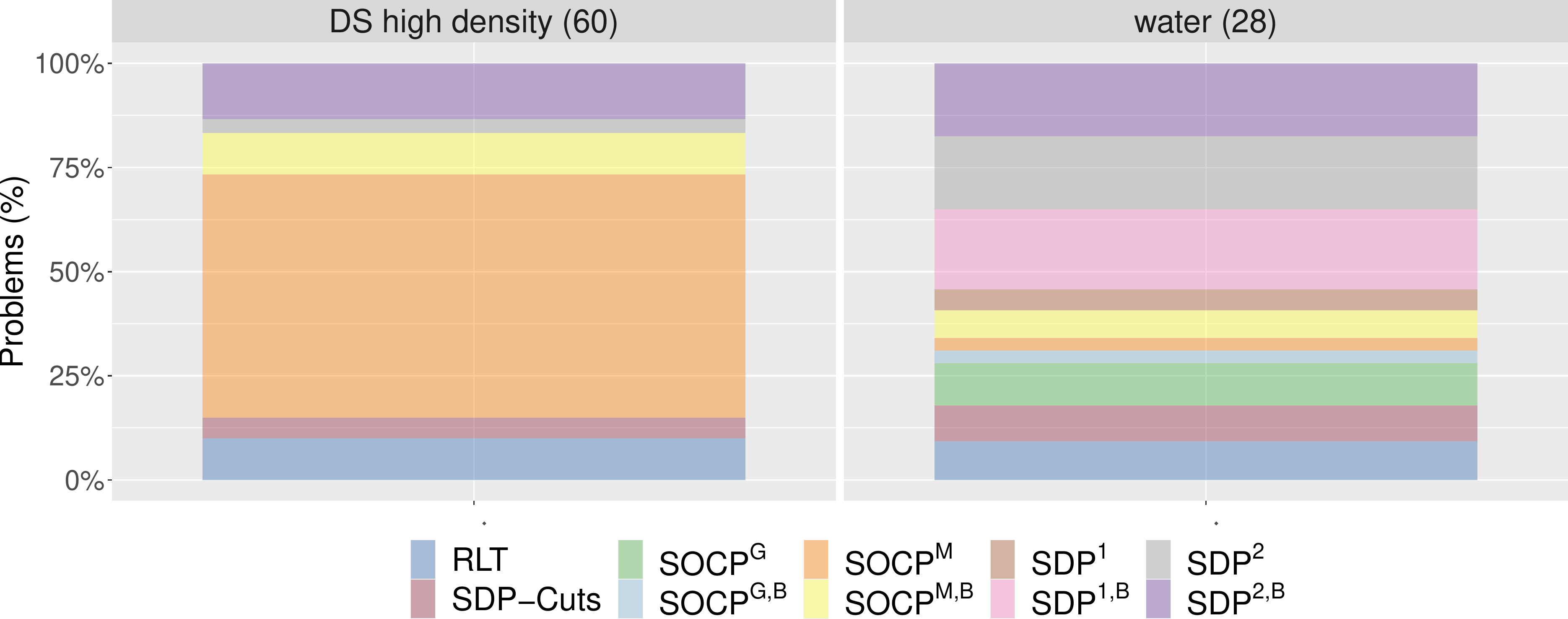}
\caption{Percentage of instances in which each version delivers the best performance in high density problems in DS and ``water''-related instances in MINLPLib.}
\label{fig:bar_comp_subsets}
\end{figure}

A natural question given the results in Figure~\ref{fig:bar_comp_all_bat} is whether or not there are specific subclasses of instances in the different test sets where a certain RLT version is noticeably dominant. A deeper analysis of the results shows that this is indeed the case, as presented in Figure~\ref{fig:bar_comp_subsets}. First, when looking at instances in DS with high density (larger that $0.5$) we can see that the versions relying on SOCP constraints and \texttt{Mosek} as a solver, \texttt{SOCP$^{M}$} and \texttt{SOCP$^{M,B}$}, are the best performing ones in around $70\%$ of the instances; \texttt{RLT} goes down to around $10\%$, \texttt{SDP-Cuts} to around $5\%$, and \texttt{SDP$^2$} and \texttt{SDP$^{2,B}$} split the remaining $15\%$. Second, we selected from MINLPLib set all instances whose name contains the word ``water'', which are problems related to the design of water networks \citep{castro2013waterdesign,teles2012globalwaterdesign} and wastewater treatment systems \citep{castro2009wastewater,castro2007wastewater}. We have that in 15 out of the 28 resulting instances, the best option is one of the following SDP-based ones: \texttt{SDP$^{1,B}$}, \texttt{SDP$^2$}, and \texttt{SDP$^{2,B}$}. Again, the instances in which \texttt{RLT} or \texttt{SDP-Cuts} are the best option are less than $20\%$. The reasons behind the notoriously good performance of the versions based on second order cone constraints for high density problems in DS and of those based on positive semidefinite constraints for ``water''-related instances in MINLPLib are definitely worth studying more deeply, but such an analysis is beyond the scope of this paper.

\begin{table}[!htbp]
    \centering
    \begin{tabular}{l|cccc}
    & \textbf{All} & \textbf{DS high density}  & \textbf{water}  \\
    \toprule\toprule
    \textbf{RLT} (across all instances) & 12.69  & 0.82 & 33994.03 \\
    \textbf{Optimal version} (instance by instance) & 6.28 & 	0.21 & 2546.88 \\
    \midrule
    \textbf{Improvement} & 50.5\% & 74.4\% & 92.5\% \\
    \bottomrule
    \end{tabular}
    \caption{Average values for $\LBpace$.}
    \label{table:resultsnoml}
\end{table}

Despite the results shown in Figure~\ref{fig:bar_comp_all_bat} and Figure~\ref{fig:bar_comp_subsets}, it is important to ensure that the variability is not spurious. For instance, it might be that all RLT versions performed very similarly to one another, which would turn most of the above discussion meaningless. In Table~\ref{table:resultsnoml} we present a concise summary of the geometric mean of $\LBpace$ for the complete set of instances and for the two special subclasses identified above. The first row measures the performance of \texttt{RLT} while the second row measures the performance of a hypothetical RLT version capable of choosing the best performing RLT version in each and every instance. The first column shows that, on aggregate on the whole set of instances, this hypothetical and optimal version would divide the pace by two, \emph{i.e.}, a substantial improvement of $50.5\%$. This improvement is much more pronounced for high density problems in DS, in which the pace gets divided by four ($74.4\%$ improvement), and even more so for ``water''-related instances in MINLPLib where the pace becomes more than ten times smaller ($92.5\%$ improvement).

\section{Machine Learning for different Conic Constraints}\label{sec:learning}

In this section we want to exploit the wide variability in the performance of the different RLT versions shown above to try to learn in advance which one should be chosen for a given instance. The goal is to design a machine learning procedure that can be trained using the different features of the instances and then choose the most promising RLT version when confronted with a new instance. The performance of the hypothetical ``Optimal version'' in Table~\ref{table:resultsnoml} represents an upper bound on the improvement that can be attained by such a machine learning version.

We follow the framework in \cite{Bissan2022}, where the authors use learning techniques to improve the performance of the RLT-based solver \texttt{RAPOSa} by learning to choose between different branching rules. The improvements reported there are substantial, with the machine learning version delivering improvements of up to $25\%$ with respect to the best original branching rule. Table~\ref{table:features} presents the full list of the input variables (features).\footnote{VIG and CMIG stand for two graphs that can be associated to any given polynomial optimization problem: \textit{variables intersection graph} and \textit{constraints-monomials intersection graph}, and whose precise definitions is given in \cite{Bissan2022}.} They capture diverse characteristics of each instance and are a key ingredient of the machine learning framework, which consists of predicting the performance ($\NLBpace$) of each branching rule on a new instance based on a regression analysis of its performance on the training instances. Then, the rule with a highest predicted performance for the given instance is chosen.

\begin{table}[!htbp]
\centering
\renewcommand{\arraystretch}{0.6}
    \begin{tabular}{ll}
      \toprule
      \multirow{5}{*}{Variables}   &  No. of variables, variance of the density of the variables \\
         & Average/median/variance of the ranges of the variables  \\
         & Average/variance of the no. of appearances of each variable\\
         & Pct. of variables not present in any monomial with degree greater than one \\
         & Pct. of variables not present in any monomial with degree greater than two \\
        \midrule
      Constraints &  No. of constraints, Pct. of equality/linear/quadratic constraints\\
    \midrule
    \multirow{3}{*}{Monomials}   &  No. of monomials\\
      & Pct. of linear/quadratic monomials, Pct. of linear/quadratic RLT variables\\
      & Average pct. of monomials in each constraint and in the objective function\\
      \midrule
      Coefficients   &  Average/variance of the coefficients\\
     \midrule
     \multirow{3}{*}{Other}   &  Degree and density of \ref{eq:PO}\\
     & No. of variables divided by no. of constrains/degree \\
     & No. of RLT variables/monomials divided by no. of constrains\\
     \midrule
    Graphs  & Density, modularity, treewidth, and transitivity of VIG and CMIG\\
    \bottomrule
    \end{tabular}
\caption{Features used for the learning.}
\label{table:features}
\end{table}

For the learning we rely on quantile regression models since as argued in \cite{Bissan2022}, the presence of outliers and of an asymmetric behavior of $\NLBpace$ (negative skewness) makes quantile regression more suitable than conventional regression models (based on the conditional mean). More precisely, we use quantile regression forests as the core tool for our analysis. Random forests were introduced in \cite{Breiman2001} as ensemble methods aggregate the information on several individual decision trees to make a single prediction. Quantile regression forests, introduced in \cite{Meinshausen2006}, are a generalization of random forests that compute an estimation of the conditional distribution of the response variable by taking into account all the observations in every leaf of every tree and not just their average. As discussed in \cite{Bissan2022}, one advantage of random forests is that the reporting of the results can be done for the complete data set in terms of Out-Of-Bag predictions; there is no need to do the reporting with respect to the underlying training and test sets. The statistical analysis is conducted in programming language~\texttt{R} \citep{Rlang}, using library \texttt{ranger}~\citep{rangerR}. Finally, the learning is conducted jointly on all sets of instances.




Table~\ref{table:resultsml} shows the remarkable improvement obtained with the machine learning (ML) version of RLT that chooses, for each given instance, the most promising version of the 10-version portfolio. 
Considering all instances, the ML version improves $32.9\%$ with respect to \texttt{RLT}, being the upper bound for learning $50.5\%$. In instances with high densities from DS test set, it improves a remarkable $69.5\%$ out of the optimal $74.4\%$, and in ``water''-related instances in MINLPLib it improves $56.7\%$ out of $92.5\%$. It is worth noting that the size of these improvements is comparable, and even superior, to those obtained in \cite{Bissan2022} when this learning scheme was introduced to learn to choose between branching rules, which shows the robustness of the proposed approach.

\begin{table}[!htbp]
    \centering
    \begin{tabular}{l|cccc}
    & \textbf{All} & \textbf{DS high density}  & \textbf{water}  \\
    \toprule\toprule
    \textbf{RLT} (across all instances) & 12.69  & 0.82 & 33994.03 \\
    \textbf{ML-based version} & 8.51 & 0.25 &	14727.22 \\
    \textbf{Optimal version} (instance by instance) & 6.28 & 	0.21 & 2546.88 \\
    \midrule
    \textbf{Improvement after learning} & 32.9\% & 69.5\% & 56.7\% \\
    \textbf{Optimal improvement} (upper bound for learning) & 50.5\% & 74.4\% & 92.5\% \\
    \bottomrule
    \end{tabular}
    \caption{Performance of the ML-based version with respect to $\LBpace$ (average across test sets).}
    \label{table:resultsml}
\end{table}

\begin{figure}[!htbp]
\centering
\includegraphics[width=0.75\textwidth]{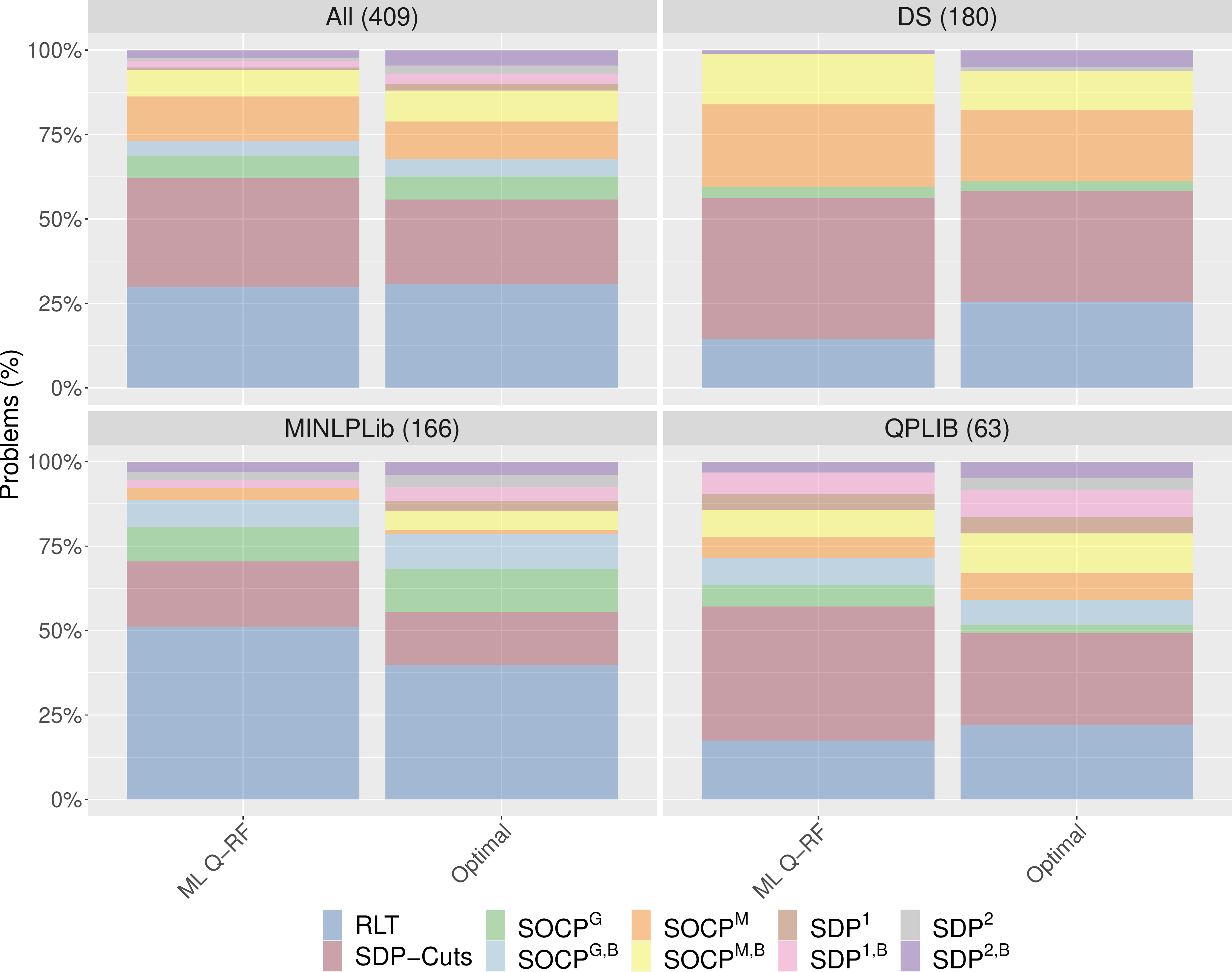}
\caption{Percentages in which each version is selected by the ML version and by the optimal one.}
\label{fig:barplot_ml_all_bat}
\end{figure}

\begin{figure}[!htbp]
\centering
\includegraphics[width=0.75\textwidth]{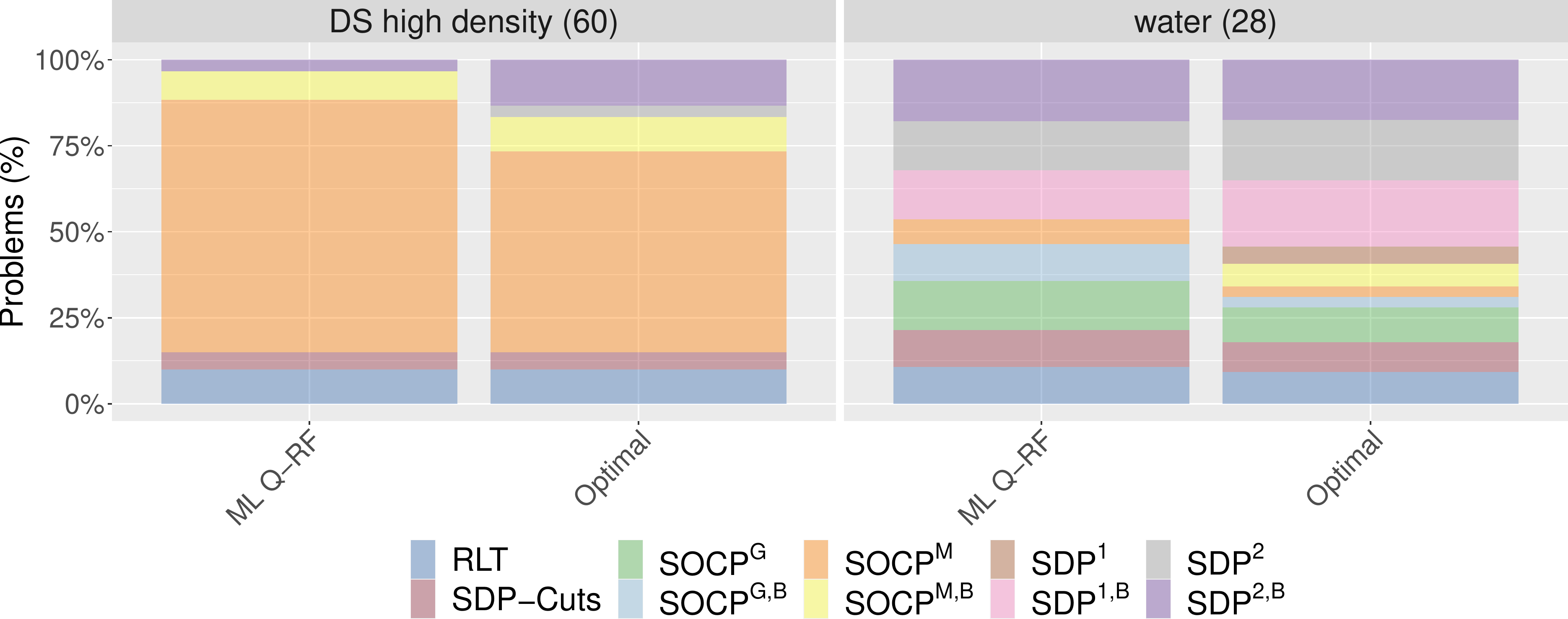}
\caption{Percentages in which each version is selected by the ML version and by the optimal one in high density problems in DS and ``water''-related instances in MINLPLib.}
\label{fig:barplot_ml_subsets}
\end{figure}

Figure~\ref{fig:barplot_ml_all_bat} represents, side by side, how often each of the ten RLT versions is selected by the ML version and by the optimal one. We can see that the behavior of the former mimics quite well that of the latter in the three sets of instances. Given that the learning is conducted jointly on the whole set of instances, the fact that the ML version adapts to the instances in DS, MINLPLib, and QPLib, is reassuring about the quality of the learning process. In particular, the SOCP versions with \texttt{Mosek} are primarily chosen for DS test set and the SDP versions are mainly chosen for QPLIB test set, the one in which they are optimal more often. In Figure~\ref{fig:barplot_ml_subsets} we further explore the behavior for the high density problems in DS and for ``water''-related instances in MINLPLib. We see that, again, the ML version mimics the patterns of the optimal version. The dominant version in DS, \texttt{SOCP$^{M}$}, is chosen in almost $75\%$ of the instances. Similarly, the three SDP dominant versions in ``water''-related instances are chosen almost $50\%$ of the time. The dominant version in DS, \texttt{SOCP$^{M}$}, is chosen in almost $75\%$ of the instances. 

\begin{figure}[!htbp]
\centering
\includegraphics[width=0.75\textwidth]{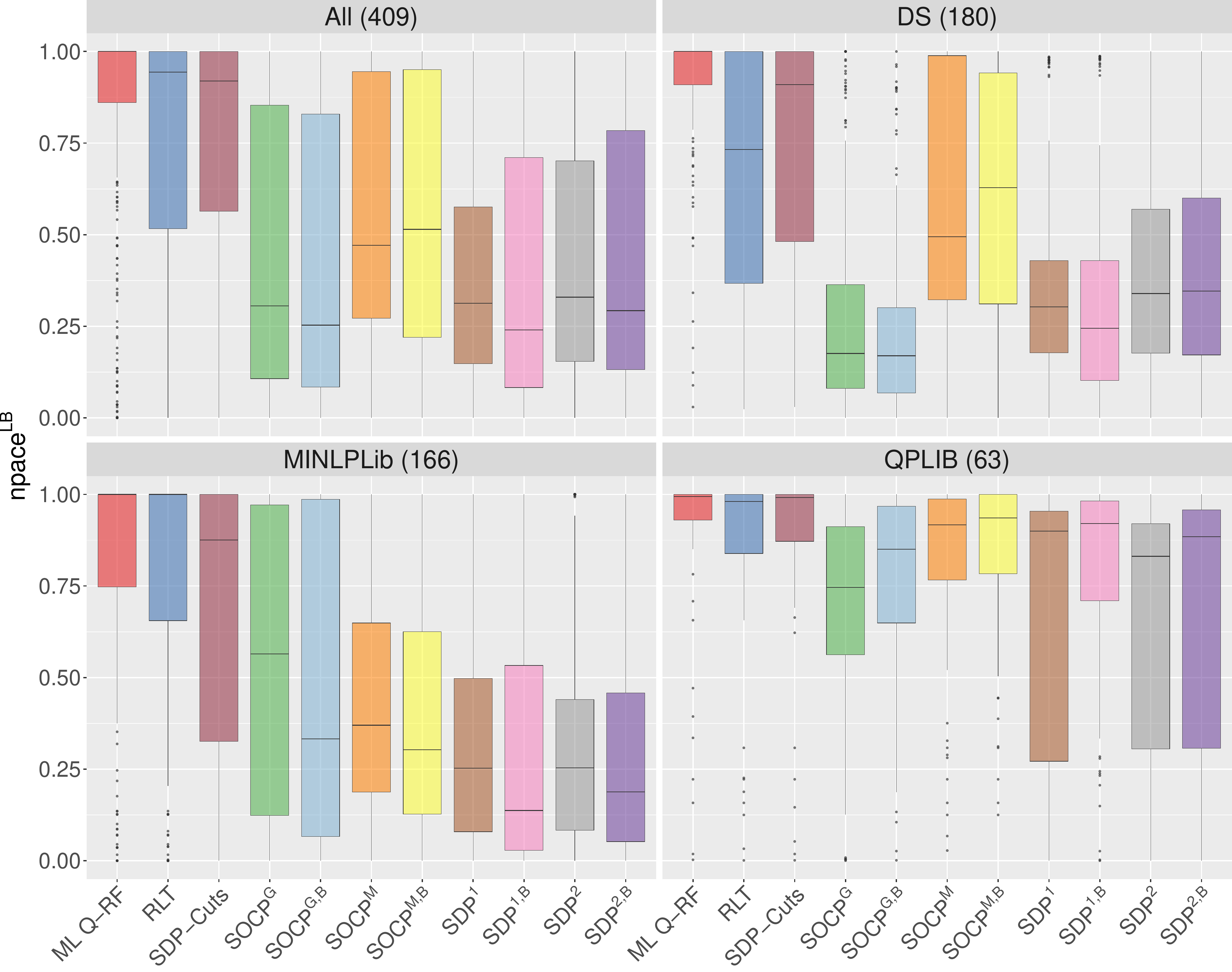}
\caption{Boxplot of $\NLBpace$ for each approach.}
\label{fig:boxplot_ml_all_bat}
\end{figure}

\begin{figure}[!htbp]
\centering
\includegraphics[width=0.75\textwidth]{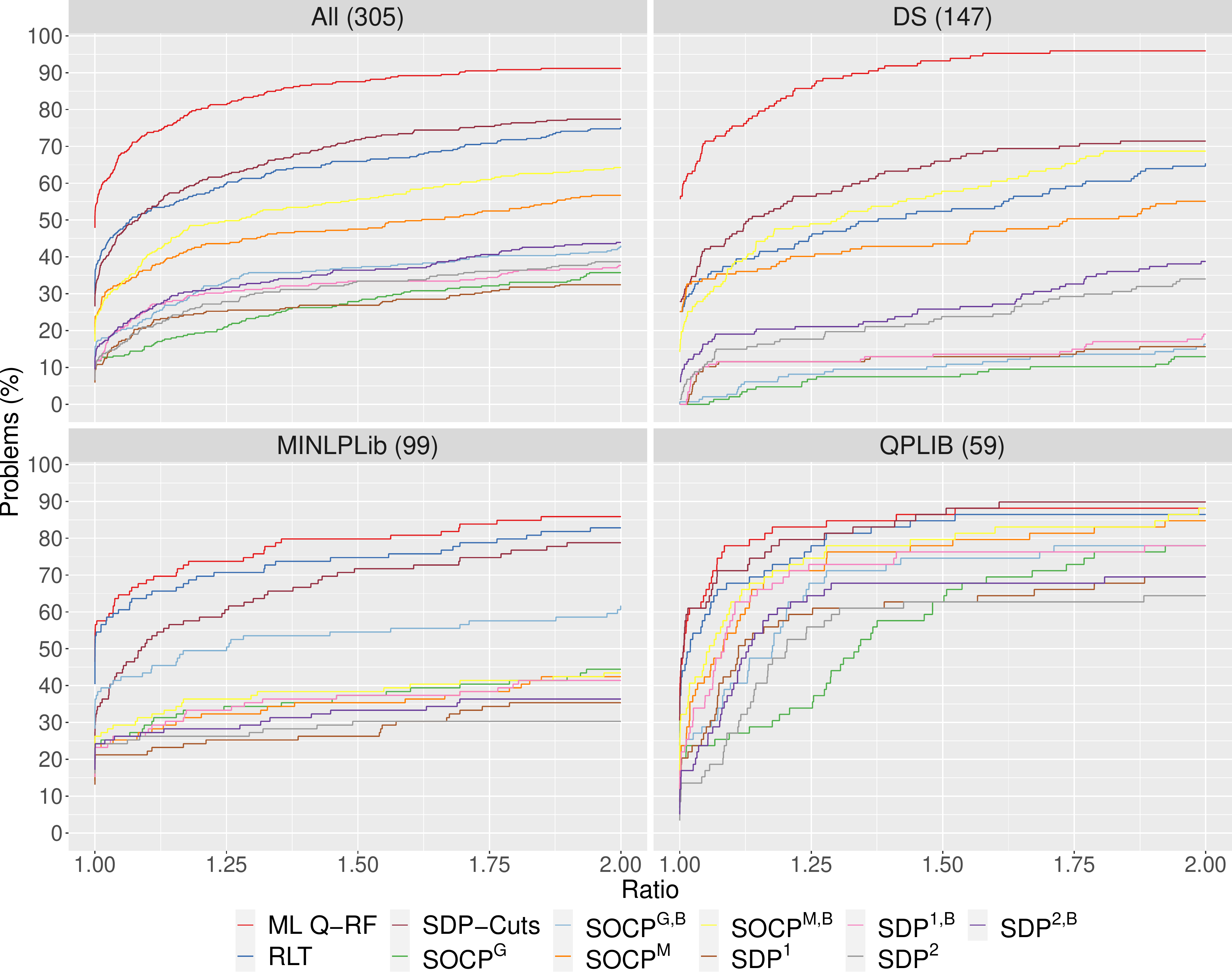}
\caption{Performance profiles of $\LBpace$ for each approach.}
\label{fig:ppml}
\end{figure}

Figure~\ref{fig:boxplot_ml_all_bat} presents boxplots summarizing the performance according to  $\NLBpace$ for all instances and for each individual test set. Recall that, by definition, values close to 1 mean that the corresponding version is almost the best one, whereas values close to 0 mean that its pace is much worse than the best one. We can see that, although \texttt{RLT} and \texttt{SDP cuts} versions are, on aggregate, the best ones in all three test sets, they are significantly outperformed by the ML version in the three of them. The improvement is particularly noticeable in DS, where the learning criterion is, by far, better than all underlying versions. This observation is further reinforced by the performance profiles (\cite{Dolan2002}) represented in Figure~\ref{fig:ppml}. Again, the ML version clearly outperforms all others, specially in DS instances.

\section{Conclusions and Future Research}\label{sec:conclusions}

The main contribution of this paper is to show that the solution to global optimality via branch-and-bound schemes of non-convex optimization problems and, in particular, polynomial optimization ones, can benefit from tightening the underlying relaxations with conic constraints. We explore different families of such constraints, building upon either second-order cones or positive semidefiniteness. We also show that the potential of these conic constraints can be successfully exploited by embedding it into a learning framework. The main goal is to predict which is the most promising type of constraints to add to the RLT relaxation at each node of the underlying branch-and-bound algorithm when confronted with a new instance. The results in Section~\ref{sec:results} show that the versions with SOCP/SDP conic constraints deliver consistently good results for instances in specific subclasses of problems: high density problems in DS and of those based on positive semidefinite constraints for ``water''-related instances in MINLPLib. 

As a future step, one may wonder to what extent one might get an even superior performance if the learning analysis was further specialized for the current setting: for example including features capturing some ``conic'' characteristics of the polynomial optimization problems and fine-tuning the regression techniques. Furthermore, an important direction for future research is to improve the understanding on the structure of these problems and the specificities that lead to the superior performance of SOCP and SDP constraints, respectively. Another direction is to investigate the number of SOCP and SDP constraints added at different nodes of the branch-and-bound tree. Additionally, we aim to extend this framework for various relaxations of polynomial optimization problems besides RLT.

\section*{Acknowledgements}
This research has been funded by FEDER and the Spanish Ministry of Science and Technology through projects MTM2014-60191-JIN, MTM2017-87197-C3, and PID2021-124030NB-C32. Brais Gonz\'alez-Rodr\'iguez acknowledges support from the Spanish Ministry of Education through FPU grant 17/02643. Raúl Alvite-Paz\'o and Samuel Alvite-Paz\'o acknowledge support from CITMAga through proyect ITMATI-R-7-JGD. Bissan Ghaddar's research is supported by Natural Sciences and Engineering Research Council of Canada Discovery Grant 2017-04185 and by the David G. Burgoyne Faculty Fellowship.

\bibliographystyle{apalike}
\bibliography{references}

\begin{thebibliography}{}

\bibitem[Andersen and Andersen, 2000]{mosek}
Andersen, E.~D. and Andersen, K.~D. (2000).
\newblock {\em The {M}osek Interior Point Optimizer for Linear Programming: An
  Implementation of the Homogeneous Algorithm}.
\newblock Springer US, Boston, MA.

\bibitem[Baltean-Lugojan et~al., 2019]{Baltean2019}
Baltean-Lugojan, R., Bonami, P., Misener, R., and Tramontani, A. (2019).
\newblock Scoring positive semidefinite cutting planes for quadratic
  optimization via trained neural networks.
\newblock Technical report, Optimization-online 7942.

\bibitem[Belotti et~al., 2009]{couenne}
Belotti, P., Lee, J., Liberti, L., Margot, F., and W{\"a}chter, A. (2009).
\newblock Branching and bounds tightening techniques for non-convex {MINLP}.
\newblock {\em Optimization Methods and Software}, 24(4-5):597--634.

\bibitem[Bengio et~al., 2021]{Bengio2021}
Bengio, Y., Lodi, A., and Prouvost, A. (2021).
\newblock Machine learning for combinatorial optimization: a methodological
  tour d'horizon.
\newblock {\em European Journal of Operational Research}, 290(2):405--421.

\bibitem[Bonami et~al., 2019]{bonami2019}
Bonami, P., Lodi, A., Schweiger, J., and Tramontani, A. (2019).
\newblock Solving quadratic programming by cutting planes.
\newblock {\em SIAM Journal on Optimization}, 29(2):1076--1105.

\bibitem[Breiman, 2001]{Breiman2001}
Breiman, L. (2001).
\newblock Random forests.
\newblock {\em Machine Learning}, 45(1):5--32.

\bibitem[Buchheim and Wiegele, 2013]{buchheim2013}
Buchheim, C. and Wiegele, A. (2013).
\newblock Semidefinite relaxations for non-convex quadratic mixed-integer
  programming.
\newblock {\em Mathematical Programming}, 141(1):435--452.

\bibitem[Burer and Vandenbussche, 2008]{burer2008}
Burer, S. and Vandenbussche, D. (2008).
\newblock A finite branch-and-bound algorithm for nonconvex quadratic
  programming via semidefinite relaxations.
\newblock {\em Mathematical Programming}, 113(2):259--282.

\bibitem[Burer and Ye, 2020]{burer2020}
Burer, S. and Ye, Y. (2020).
\newblock Exact semidefinite formulations for a class of (random and
  non-random) nonconvex quadratic programs.
\newblock {\em Mathematical Programming}, 181(1):1--17.

\bibitem[Bussieck et~al., 2003]{minlplib}
Bussieck, M.~R., Drud, A.~S., and Meeraus, A. (2003).
\newblock {MINLPL}ib-a collection of test models for mixed-integer nonlinear
  programming.
\newblock {\em INFORMS Journal on Computing}, 15:114--119.

\bibitem[Castro et~al., 2007]{castro2007wastewater}
Castro, P.~M., Matos, H.~A., and Novais, A.~Q. (2007).
\newblock An efficient heuristic procedure for the optimal design of wastewater
  treatment systems.
\newblock {\em Resources, conservation and recycling}, 50(2):158--185.

\bibitem[Castro and Teles, 2013]{castro2013waterdesign}
Castro, P.~M. and Teles, J.~P. (2013).
\newblock Comparison of global optimization algorithms for the design of
  water-using networks.
\newblock {\em Computers \& chemical engineering}, 52:249--261.

\bibitem[Castro et~al., 2009]{castro2009wastewater}
Castro, P.~M., Teles, J.~P., and Novais, A.~Q. (2009).
\newblock Linear program-based algorithm for the optimal design of wastewater
  treatment systems.
\newblock {\em Clean Technologies and Environmental Policy}, 11(1):83--93.

\bibitem[Dalkiran and Sherali, 2013]{Sherali2013}
Dalkiran, E. and Sherali, H.~D. (2013).
\newblock Theoretical filtering of {RLT} bound-factor constraints for solving
  polynomial programming problems to global optimality.
\newblock {\em Journal of Global Optimization}, 57(4):1147--1172.

\bibitem[Dalkiran and Sherali, 2016]{Dalkiran2016}
Dalkiran, E. and Sherali, H.~D. (2016).
\newblock {RLT}-{POS}: Reformulation-linearization technique-based optimization
  software for solving polynomial programming problems.
\newblock {\em Mathematical Programming Computation}, 8:337--375.

\bibitem[Dolan and Mor\'e, 2002]{Dolan2002}
Dolan, E.~D. and Mor\'e, J.~J. (2002).
\newblock Benchmarking optimization software with performance profiles.
\newblock {\em Mathematical Programming}, 91:201--213.

\bibitem[Elloumi and Lambert, 2019]{elloumi2019}
Elloumi, S. and Lambert, A. (2019).
\newblock Global solution of non-convex quadratically constrained quadratic
  programs.
\newblock {\em Optimization methods and software}, 34(1):98--114.

\bibitem[{FICO}, 2022]{xpress}
{FICO} (2022).
\newblock {FICO} {X}press {O}ptimization {S}uite.
\newblock Available at:
  \url{https://www.fico.com/en/products/fico-xpress-optimization}.

\bibitem[Furini et~al., 2018]{qplib}
Furini, F., Traversi, E., Belotti, P., Frangioni, A., Gleixner, A., Gould, N.,
  Liberti, L., Lodi, A., Misener, R., Mittelmann, H., Sahinidis, N., Vigerske,
  S., and Wiegele, A. (2018).
\newblock {QPLIB}: a library of quadratic programming instances.
\newblock {\em Mathematical Programming Computation}, 1:237--265.

\bibitem[Ghaddar et~al., 2011a]{ghaddar2011max}
Ghaddar, B., Anjos, M.~F., and Liers, F. (2011a).
\newblock A branch-and-cut algorithm based on semidefinite programming for the
  minimum k-partition problem.
\newblock {\em Annals of Operations Research}, 188(1):155--174.

\bibitem[Ghaddar et~al., 2022]{Bissan2022}
Ghaddar, B., Gómez-Casares, I., González-Díaz, J., González-Rodríguez, B.,
  Pateiro-López, B., and Rodríguez-Ballesteros, S. (2022).
\newblock Learning for spatial branching: An algorithm selection approach.
\newblock Technical report.

\bibitem[Ghaddar and Jabr, 2019]{ghaddar2019}
Ghaddar, B. and Jabr, R.~A. (2019).
\newblock Power transmission network expansion planning: A semidefinite
  programming branch-and-bound approach.
\newblock {\em European Journal of Operational Research}, 274(3):837--844.

\bibitem[Ghaddar et~al., 2011b]{ghaddar2011}
Ghaddar, B., Vera, J.~C., and Anjos, M.~F. (2011b).
\newblock Second-order cone relaxations for binary quadratic polynomial
  programs.
\newblock {\em SIAM Journal on Optimization}, 21(1):391--414.

\bibitem[Gonz{\'a}lez-Rodr{\'\i}guez et~al., 2020]{Gonzalez-rodriguez2020}
Gonz{\'a}lez-Rodr{\'\i}guez, B., Ossorio-Castillo, J., Gonz{\'a}lez-D{\'\i}az,
  J., Gonz{\'a}lez-Rueda, {\'A}.~M., Penas, D.~R., and
  Rodr{\'\i}guez-Mart{\'\i}nez, D. (2020).
\newblock Computational advances in polynomial optimization: {RAPOSa}, a freely
  available global solver.
\newblock Technical report, Optimization-online 7942.

\bibitem[{Gurobi Optimization}, 2022]{gurobi}
{Gurobi Optimization} (2022).
\newblock Gurobi {O}ptimizer {R}eference {M}anual.
\newblock Available at: \url{http://www.gurobi.com}.

\bibitem[{IBM Corp.}, 2022]{cplex}
{IBM Corp.} (2022).
\newblock {IBM} {ILOG} {CPLEX} {O}ptimization {S}tudio. {CPLEX} {U}ser’s
  {M}anual.
\newblock Available at:
  \url{https://www.ibm.com/es-es/products/ilog-cplex-optimization-studio}.

\bibitem[Krislock et~al., 2017]{krislock2017}
Krislock, N., Malick, J., and Roupin, F. (2017).
\newblock Biqcrunch: A semidefinite branch-and-bound method for solving binary
  quadratic problems.
\newblock {\em ACM Trans. Math. Softw.}, 43(4).

\bibitem[Lasserre, 2001]{lasserre2001}
Lasserre, J.~B. (2001).
\newblock Global optimization with polynomials and the problem of moments.
\newblock {\em SIAM Journal on optimization}, 11(3):796--817.

\bibitem[Lodi and Zarpellon, 2017]{Lodi2017}
Lodi, A. and Zarpellon, G. (2017).
\newblock On learning and branching: a survey.
\newblock {\em Top}, 25(2):207--236.

\bibitem[Meinshausen, 2006]{Meinshausen2006}
Meinshausen, N. (2006).
\newblock Quantile regression forests.
\newblock {\em Journal of Machine Learning Research}, 7:983--999.

\bibitem[{MOSEK ApS}, 2022]{mosekmanual}
{MOSEK ApS} (2022).
\newblock {\em Introducing the {MOSEK} {O}ptimization {S}uite 9.3.20}.

\bibitem[Parrilo, 2003]{parrilo2003}
Parrilo, P.~A. (2003).
\newblock Semidefinite programming relaxations for semialgebraic problems.
\newblock {\em Mathematical programming}, 96(2):293--320.

\bibitem[Piccialli et~al., 2022]{piccialli2022}
Piccialli, V., Sudoso, A.~M., and Wiegele, A. (2022).
\newblock Sos-sdp: an exact solver for minimum sum-of-squares clustering.
\newblock {\em INFORMS Journal on Computing}.

\bibitem[{R Core Team}, 2021]{Rlang}
{R Core Team} (2021).
\newblock {\em R: A Language and Environment for Statistical Computing}.
\newblock R Foundation for Statistical Computing, Vienna, Austria.

\bibitem[Rendl et~al., 2010]{rendl2010}
Rendl, F., Rinaldi, G., and Wiegele, A. (2010).
\newblock Solving max-cut to optimality by intersecting semidefinite and
  polyhedral relaxations.
\newblock {\em Mathematical Programming}, 121(2):307--335.

\bibitem[Sherali et~al., 2012]{Sherali2012cuts}
Sherali, H.~D., Dalkiran, E., and Desai, J. (2012).
\newblock Enhancing {RLT}-based relaxations for polynomial programming problems
  via a new class of $v$-semidefinite cuts.
\newblock {\em Computational Optimization and Applications}, 52(2):483–506.

\bibitem[Sherali and Tuncbilek, 1992]{Sherali1992}
Sherali, H.~D. and Tuncbilek, C.~H. (1992).
\newblock A global optimization algorithm for polynomial programming problems
  using a reformulation-linearization technique.
\newblock {\em Journal of Global Optimization}, 2(1):101--112.

\bibitem[Shor, 1987]{shor1987}
Shor, N.~Z. (1987).
\newblock An approach to obtaining global extremums in polynomial mathematical
  programming problems.
\newblock {\em Cybernetics}, 23(5):695--700.

\bibitem[Teles et~al., 2012]{teles2012globalwaterdesign}
Teles, J.~P., Castro, P.~M., and Matos, H.~A. (2012).
\newblock Global optimization of water networks design using multiparametric
  disaggregation.
\newblock {\em Computers \& Chemical Engineering}, 40:132--147.

\bibitem[Wright and Ziegler, 2017]{rangerR}
Wright, M.~N. and Ziegler, A. (2017).
\newblock {ranger}: A fast implementation of random forests for high
  dimensional data in {C++} and {R}.
\newblock {\em Journal of Statistical Software}, 77(1):1--17.

\end{thebibliography}
\end{document}